\documentclass{ifacconf}

\usepackage{graphicx}      
\usepackage{natbib}        
\usepackage{epsfig} 
\usepackage{amsmath} 
\usepackage{amssymb} 
\usepackage[algoruled,vlined,algo2e]{algorithm2e}
\usepackage{multicol}
\usepackage{subcaption}
\usepackage{algorithm}
\usepackage{algorithmic}
\usepackage{verbatimbox}

\usepackage{booktabs}

\begin{document}
\begin{frontmatter}

\title{Nonlinear Model Predictive Control\\and System Identification\\ for a Dual-hormone Artificial Pancreas} 

\author{Asbj{\o}rn Thode Reenberg\textsuperscript{*},}
\author{Tobias K. S. Ritschel\textsuperscript{*},}
\author{Emilie B. Lindkvist\textsuperscript{**},}
\author{Christian Laugesen\textsuperscript{**},}
\author{Jannet Svensson\textsuperscript{**},}
\author{Ajenthen G. Ranjan\textsuperscript{**},}
\author{Kirsten N{\o}rgaard\textsuperscript{**}},
\author{John Bagterp J{\o}rgensen\textsuperscript{*}}
	\thanks{
		Corresponding author: J. B. J{\o}rgensen (E-mail: {\tt\small jbjo@dtu.dk}).}

\address[info]{Department of Applied Mathematics and Computer Science, Technical University of Denmark, DK-2800 Kgs. Lyngby, Denmark.}
\address[info]{Steno Diabetes Center Copenhagen, Clinical Research, DK-2730 Herlev, Denmark.}

\begin{abstract}                
In this work, we present a switching nonlinear model predictive control (NMPC) algorithm for a dual-hormone artificial pancreas (AP), and we use maximum likelihood estimation (MLE) to identify model parameters. A dual-hormone AP consists of a continuous glucose monitor (CGM), a control algorithm, an insulin pump, and a glucagon pump. The AP is designed with a heuristic to switch between insulin and glucagon as well as state-dependent constraints. We extend an existing glucoregulatory model with glucagon and exercise for simulation, and we use a simpler model for control. We test the AP (NMPC and MLE) using in silico numerical simulations on 50 virtual people with type 1 diabetes. The system is identified for each virtual person based on data generated with the simulation model. The simulations show a mean of 89.3\% time in range (3.9--10~mmol/L) and no hypoglycemic events.

\end{abstract}


\end{frontmatter}

\section{Introduction}
Type 1 diabetes (T1D) is a chronic metabolic disorder which prevents the pancreas from producing insulin. People with T1D require life-long treatment with daily injections of insulin in order to prevent hyperglycemia (i.e., high blood glucose concentrations). Prolonged hyperglycemia leads to a range of health complications, e.g., cardiovascular disease, chronic kidney disease, and damage to the nerves and eyes. Over 9\% of the world population suffers from diabetes, 5-10\% of those have T1D, and 10\% of the 2019 global health expenditure (USD~760~billion) was spent on diabetes~\citep{IDF:2019}.

The treatment of T1D is tedious and time-consuming, and if managed poorly, it can lead to both hyper- and hypoglycemia (low blood glucose concentrations). Hypoglycemia can, in severe cases, cause a variety of acute complications including loss of consciousness, seizures, and death.
Therefore, there is a significant interest in developing closed-loop diabetes treatment systems based on feedback control. Such systems are referred to as \emph{artificial pancreases}~(APs). They consist of 1)~a continuous glucose monitor (CGM) (the sensor), 2)~a control algorithm (e.g., implemented on a smartphone or a dedicated device), and 3)~an insulin pump (the actuator).
Many control strategies have been proposed for this purpose, including fuzzy logic~\citep{Biester:etal:2019}, proportional-integral-derivative~(PID) control~\citep{Sejersen:etal:2021, Huyett:etal:2015, Jorgensen:etal:2019}, and model predictive control (MPC). MPC is a closed-loop feedback strategy that uses the moving horizon optimization principle, i.e., it involves solving a sequence of open-loop optimal control problems (OCPs). Both linear MPC~(LMPC)~\citep{Chakrabarty:etal:2020, Messori:etal:2018}, and nonlinear MPC~(NMPC)~\citep{Hovorka:etal:2004, Boiroux1:etal:2018, Boiroux2:etal:2018} have been considered.

Most algorithms are designed for \emph{single-hormone} systems where only insulin is administered. Consequently, they are unable to actively counteract low blood sugar concentrations which can occur in a variety of situations, e.g., in connection with physical activity. Therefore, researchers currently investigate \emph{dual-hormone} systems that administer both insulin and glucagon~\citep{Peters:Haidar:2018, Infante:etal:2021}. In contrast to insulin, glucagon causes an increase in the blood glucose level. \citet{Moscardo:etal:2019} develop a dual-hormone control algorithm based on proportional-derivative~(PD) control, and \citet{Boiroux3:etal:2018} develop LMPC algorithms based on a variety of different transfer function models. Other hormones than glucagon have also been considered, e.g., pramlintide which slows down gastric emptying~\citep{Haidar:etal:2020}.

In this work, we present a dual-hormone NMPC algorithm for administering insulin and glucagon. The algorithm is based on an extended Medtronic Virtual patient (MVP) model~\citep{Kanderian:etal:2009}, and we use maximum likelihood estimation (MLE) to identify the model parameters. We use the continuous-discrete extended Kalman filter (CD-EKF) in both the parameter estimation and in the NMPC algorithm. Furthermore, we use a state-dependent heuristic for switching between administering insulin and glucagon (which cannot be administered simultaneously).
We use an extension of the model developed by \citet{Hovorka:etal:2002} to perform closed-loop simulations. It is extended with 1)~a model of the measurement delay of the CGM~\citep{Facchinetti:etal:2014}, 2)~a pharmacokinetic model of subcutaneous glucagon injection~\citep{Haidar:etal:2013}, and 3)~a model of the effect of physical activity~\citep{Rashid:etal:2019}. We present numerical results for 50~virtual people with T1D and demonstrate that the AP, including the parameter estimation, satisfies the time in range (TIR) targets described by \citet{Holt:etal:2021b}.

The remainder of the paper is structured as follows. We present the extension of the model by Hovorka et al. in Section~\ref{sec:SimulationModel} and the MVP model in Section~\ref{sec:ControlModel}. In Section~\ref{sec:ParameterEstimation}, we describe the parameter estimation problem, and we describe the NMPC algorithm in Section~\ref{sec:OptimalControlProblem}. We present the state-dependent switching heuristic in Section~\ref{sec:Heuristics}, and we discuss the numerical results in Section~\ref{sec:Results}. Finally, conclusions are given in Section~\ref{sec:Conclusion}.

\section{Simulation Model}
\label{sec:SimulationModel}
The model by \citet{Hovorka:etal:2002, Hovorka:etal:2004} consists of an insulin subsystem, a meal subsystem, and a glucose subsystem. We extend it with a glucagon subsystem~\citep{Haidar:etal:2013}, an exercise subsystem~\citep{Rashid:etal:2019}, and a CGM subsystem~\citep{Facchinetti:etal:2014}.
\subsection{Insulin subsystem}
The insulin absorption and insulin concentration are described by
\begin{subequations}\label{eq:hovorka:insulin}
	\begin{align}
		\label{eq:hovorka:insulin:s1}
		\dot S_1(t) &= u_I(t)    - \frac{S_1(t)}{\tau_S}, \\
		\label{eq:hovorka:insulin:s2}
		\dot S_2(t) &= \frac{S_1(t)}{\tau_S} - \frac{S_2(t)}{\tau_S}, \\
		\label{eq:hovorka:insulin:i}
		\dot I(t)   &= \frac{1}{V_I}\frac{S_2(t)}{\tau_S} - k_e I(t),
	\end{align}
\end{subequations}
where $S_1$ and $S_2$~[mU] are a two compartment chain representing the insulin absorption, $\tau_S$~[min] is the insulin absorption time constant, $I$~[mU/L] is the plasma insulin concentration, $V_I$~[L] is the insulin distribution volume, and $k_e$~[1/min] is the elimination rate. $u_I = u_{ba} + u_{bo}$ [mU/min] is the insulin infusion rate where $u_{ba}$ [mU/min] is the basal infusion rate and $u_{bo}$ [mU/min] is the bolus infusion rate.

\subsection{Insulin action subsystem}
The insulin action on the glucose kinetics is described by
\begin{subequations}\label{eq:hovorka:insulin:action}
	\begin{align}
		\label{eq:hovorka:insulin:action:x1}
		\dot x_1(t) &= k_{b1} I(t) - k_{a1} x_1(t), \\
		\label{eq:hovorka:insulin:action:x2}
		\dot x_2(t) &= k_{b2} I(t) - k_{a2} x_2(t), \\
		\label{eq:hovorka:insulin:action:x3}
		\dot x_3(t) &= k_{b3} I(t) - k_{a3} x_3(t),
	\end{align}
\end{subequations}
where $x_1$~[1/min], $x_2$~[1/min], $x_3$~[1/min] represent the effects of insulin on the glucose distribution, the glucose disposal, and the endogenous glucose production. Furthermore, $k_{bi}$~[(L/mU)/min$^2$] and $k_{ai}$~[1/min] for $i = 1, 2, 3$ are the activation and deactivation rate constants.

\subsection{Meal subsystem}
The meal absorption subsystem is represented by 
\begin{subequations}\label{eq:hovorka:meal}
	\begin{align}
		\label{eq:hovorka:meal:d1}
		\dot D_1(t) &= A_G D(t)    - \frac{D_1(t)}{\tau_D}, \\
		\label{eq:hovorka:meal:d2}
		\dot D_2(t) &= \frac{D_1(t)}{\tau_D} - \frac{D_2(t)}{\tau_D},
	\end{align}
\end{subequations}
where $D_1$~[mmol] and $D_2$~[mmol] are a two compartment chain representing the meal absorption, $D$~[mmol/min] is the meal carbohydrate content, $A_G$~[--] is the carbohydrate bioavailability, and $\tau_D$~[min] is the meal time constant.

\subsection{Glucagon subsystem}
The glucagon subsystem is described by
\begin{subequations}\label{eq:hovorka:glucagon}
	\begin{align}
		\label{eq:hovorka:glucagon:q1g}
		\dot Q_1^G(t) &= u_G(t) - \frac{Q_1^G(t)}{\tau_{Glu}}, \\
		\label{eq:hovorka:glucagon:q2g}
		\dot Q_2^G(t) &= \frac{Q_1^G(t)}{\tau_{Glu}} - \frac{Q_2^G(t)}{\tau_{Glu}},
	\end{align}
\end{subequations}
where $Q_1^G$~[$\mu$g] and $Q_2^G$~[$\mu$g] are a two-compartment chain representing the glucagon absorption, $u_G$~[$\mu$g/min] is the glucagon infusion rate, and $\tau_{Glu}$~[min] is a time constant.

\subsection{Exercise subsystem}
The exercise model describes the increased glucose consumption and insulin sensitivity during and after physical activity. The effect of high-intensity exercise is not modeled. The exercise subsystem consists of
\begin{subequations}\label{eq:hovorka:exercise}
	\begin{align}
		\label{eq:hovorka:exercise:e1}
		\dot E_1(t) &= \frac{HR(t) - HR_0 - E_1(t)}{\tau_{HR}}, \\
		\label{eq:hovorka:exercise:te}
		\dot T_E(t) &= \frac{c_1 f_{E1}(t) + c_2 - T_E(t)}{\tau_{ex}}, \\
		\label{eq:hovorka:exercise:e2}
		\dot E_2(t) &= -\left(\frac{f_{E1}(t)}{\tau_{in}} + \frac{1}{T_E(t)}\right) E_2(t) + \frac{f_{E1}(t) T_E(t)}{c_1 + c_2}, \\
		\label{eq:hovorka:exercise:fe1}
		f_{E1}(t) &= \frac{\left(\frac{E_1(t)}{a\cdot HR_0}\right)^n}{1 + \left(\frac{E_1(t)}{a\cdot HR_0}\right)^n},
	\end{align}
\end{subequations}
where $E_1$~[BPM] is the short-term effect, $E_2$~[min] is the long-term effect, $T_E$ [min] is the characteristic time for the long-term effect, $HR$~[BPM] is the heart rate, $HR_0$~[BPM] is the resting heart rate, $\tau_{HR}$~[min] is the time constant, $c_1$~[min] and $c_2$~[min] define the steady state value for $T_E$, $\tau_{ex}$~[min] is the time constant for how fast $T_E$ reaches steady state, and $a$~[--], $n$~[--], and $\tau_{in}$~[min] specify the intensity and time constant of the long-term effect on the insulin sensitivity.

\subsection{Glucose subsystem}
The glucose kinetics are represented by
\begin{subequations}\label{eq:hovorka:glucose}
	\begin{align}
		\label{eq:hovorka:glucose:q1}
		\dot Q_1(t) &= \frac{D_2(t)}{\tau_D} - F_{01,c}(t) - F_R(t) - x_1(t) Q_1(t) \\ \nonumber &+ k_{12} Q_2(t) + EGP(t) + Q_G(t) - Q_{E21}(t), \\
		\label{eq:hovorka:glucose:q2}
		\dot Q_2(t) &= x_1(t) Q_1(t) - k_{12} Q_2(t) - x_2 Q_2(t) \\ \nonumber & +Q_{E21}(t)-Q_{E22}(t) - Q_{E1}(t),
	\end{align}
\end{subequations}
where
\begin{subequations}\label{eq:hovorka:glucose:auxiliary}
	\begin{align}
		\label{eq:hovorka:glucose:auxiliary:egp}
		EGP(t) &= EGP_0 (1 - x_3(t)), \\
		\label{eq:hovorka:glucose:auxiliary:qg}
		Q_G(t) &= K_{Glu} V_G Q_2^G(t), \\
		\label{eq:hovorka:glucose:auxiliary:qe21}
		Q_{E21}(t) &= \alpha E_2(t)^2 x_1(t) Q_1(t), \\
		\label{eq:hovorka:glucose:auxiliary:qe22}
		Q_{E22}(t) &= \alpha E_2(t)^2 x_2(t) Q_2(t), \\
		\label{eq:hovorka:glucose:auxiliary:qe1}
		Q_{E1}(t)  &= \beta \frac{E_1(t)}{HR_0},
	\end{align}
\end{subequations}
and
\begin{subequations}\label{eq:hovorka:glucose:piecewise}
	\begin{align}
		\label{eq:hovorka:glucose:piecewise:f01c}
		F_{01,c}(t) &=
		\begin{cases}
			F_{01} & G(t) \geq 4.5~\mathrm{mmol/L}, \\
			F_{01} G(t)/4.5 & \text{otherwise},
		\end{cases} \\
		\label{eq:hovorka:glucose:piecewise:fr}
		F_R(t) &=
		\begin{cases}
			0.003(G(t) - 9) V_G & G(t) \geq 9~\mathrm{mmol/L}, \\
			0 & \text{otherwise},
		\end{cases} \\
		G(t) &= \frac{Q_1(t)}{V_G},
	\end{align}
\end{subequations}
where $Q_1$~[mmol] and $Q_2$~[mmol] represent the accessible and non-accessible compartments, $k_{12}$~[1/min] is the transfer rate, $F_{01}$~[mmol/min] and $F_{01,c}$~[mmol/min] are the nominal and corrected total non-insulin dependent glucose flux, $F_R$~[mmol/min] is the renal glucose clearance,  $V_G$~[L] is the glucose distribution volume, $EGP_0$~[mmol] is the endogenous glucose production, $K_{Glu}$~[(mmol/L)/$\mu$g/min] is the glucagon gain, $\alpha$~[1/min$^2$] is the exercise-induced insulin action, $\beta$~[mmol/min] is the exercise-induced insulin-independent glucose uptake rate, and $G$~[mmol/L] is the glucose concentration.

\subsection{CGM subsystem}
The CGM subsystem describes the glucose transfer from the plasma to the interstitial tissue by
\begin{equation}\label{eq:hovorka:cgm}
	\dot G_I(t) = \frac{G(t)}{\tau_{IG}} - \frac{G_I(t)}{\tau_{IG}}, 
\end{equation}
where $G_I$~[mmol/L] is the interstitial glucose concentration and $\tau_{IG}$ is a time constant.

\section{Control model}
\label{sec:ControlModel}
In the AP, we use an extension of the Medtronic Virtual patient (MVP) model~\citep{Kanderian:etal:2009} represented as a system of coupled stochastic differential equations. The MVP model describes the glucose-insulin dynamics and we extend it with the meal subsystem, the glucagon subsystem and the CGM subsystem from the simulation model described in Section \ref{sec:SimulationModel}. 

\subsection{Insulin subsystem}
The insulin absorption subsystem consists of
\begin{subequations}\label{eq:mvp:insulin}
	\begin{align}
		\label{eq:mvp:insulin:isc}
		dI_{SC}(t) 	&= k_1 \left(\frac{u_{ba}(t) + u_{bo}(t)}{C_I} - I_{SC}(t)\right)dt, \\
		\label{eq:mvp:insulin:ip}
		dI_P(t) 		&= k_2 \left(I_{SC}(t) - I_P(t)\right)dt,
	\end{align}
\end{subequations}
where $I_{SC}$~[mU/L] is the subcutaneous insulin concentration, $I_P$~[mU/L] is the plasma insulin concentration, $k_2=k_1$~[1/min] is the inverse insulin absorption time constant, and $C_I$~[L/min] is the insulin clearance rate.

		%
		%
		%

\subsection{Glucose subsystem}
Here, we describe the insulin effect, the blood glucose concentration and insulin sensitivity. The blood glucose concentration and insulin sensitivity is modeled as stochastic differential equations:
\begin{subequations}\label{eq:mvp:glucose}
	\begin{align}
		\label{eq:mvp:glucose:ieff}
		dI_{EFF}(t) 	&= p_2 \left(S_I(t) I_P(t) - I_{EFF}(t)\right)dt, \\
		\label{eq:mvp:glucose:g}
		dG(t)  		&= [-(GEZI + I_{EFF}(t)) G(t) + EGP\\ \nonumber
		&+ R_A(t) + K_{Glu} Q_2^G(t)]dt + \sigma_G dw_G(t),\\
		\label{eq:mvp:glucose:logsi}
		d\log(S_I(t)) &= \sigma_{S_I} dw_{S_I}(t),
	\end{align}
\end{subequations}
where $I_{EFF}$~[1/min] is the insulin effect, $p_2 = k_1$~[1/min] is the inverse insulin action time constant, $S_I$~[(L/mU)/min] is the insulin sensitivity, $GEZI$~[1/min] is the glucose effectiveness, $EGP$~[(mmol/L)/min] is the endogenous glucose production, $\sigma_G$ and $\sigma_{S_I}$ are the glucose and insulin sensitivity diffusion coefficients, and $w_G$ and $w_{S_I}$ are standard Wiener processes.
The meal rate of appearance, $R_A$~[(mmol/L)/min], is
\begin{equation}
	R_A(t) = \frac{k_m D_2(t)}{V_G},
\end{equation}
where $k_m$ [1/min] is a time constant.




\section{Parameter Estimation}
\label{sec:ParameterEstimation}
We use MLE based on the CD-EKF to estimate the parameters in the MVP model
given $N+1$ CGM measurements of the blood glucose concentration, $\mathcal Y_N = \{y_0,y_1,\dots, y_N\}$.
The MVP model is in the form
\begin{subequations}\label{eq:mle:system}
	\begin{align}
		\label{eq:mle:system:dx}
		dx(t) &= f(t, x(t), u(t), d(t), \theta) dt + \sigma(\theta) dw(t),\\
		\label{eq:mle:system:y}
		y_k   &= g(t_k, x_k, \theta) + v_k,
	\end{align}
\end{subequations}
where $t$ is the time, $x$ are the states, $u$ are the manipulated inputs, $d$ are the disturbances, $\theta$ are the parameters, and $y_k = y(t_k)$ are the measured variables. $\sigma$ is the diffusion coefficient, $w$ is a standard Wiener process (i.e., $dw(t) \sim N(0, \mathrm{I} dt)$), and $v_k \sim N(0, R)$ is the measurement noise.

The MLE of the parameters, $\hat \theta$, is given by
\begin{equation}\label{eq:mle:thetahat}
	\hat \theta = \arg \min_{\theta} V(\theta),
\end{equation}
where $V$ is the negative log-likelihood function:
\begin{equation}\label{eq:mle:optimization:problem}
	V(\theta) = -\log p(\mathcal{Y}_N | \theta).
\end{equation}
Here, $p(\mathcal{Y}_N | \theta)$ is the conditional probability density function of the stochastic observations in the system~\eqref{eq:mle:system} evaluated at the observed blood glucose concentrations for a given set of parameters, $\theta$.
The negative log-likelihood function is given by
\begin{align}\label{eq:mle:v}
	V(\theta)
	&= \frac{(N+1) n_y}{2} \log(2\pi) + \frac{1}{2}\sum_{k=0}^{N}\log[\det(R_{e,k}(\theta))] \nonumber \\
	&+ e_k(\theta)^T [R_{e,k}(\theta)]^{-1} e_k(\theta),
\end{align}
where $n_y = 1$ and $e_k(\theta)$ is the innovation:
\begin{equation}\label{eq:mle:ek}
	e_k(\theta) = y_k - \hat{y}_{k|k-1}(\theta).
\end{equation}
Given an initial estimate of the states (which is also estimated) and their covariance, $\hat x_{0}(\theta)$ and $P_0$, we use the CD-EKF to compute the one-step predictions of the observed variables, $\hat{y}_{k|k-1}(\theta)$, and the covariance of the innovations, $R_{e, k}(\theta)$. We refer to the paper by~\citet{Boiroux:etal:2019} for more details.
%
%
%
%

\section{Nonlinear model predictive control}
\label{sec:OptimalControlProblem}
The NMPC algorithm receives a CGM measurement of the blood glucose concentration every 5 minutes. Subsequently, the CD-EKF is used to compute a filtered estimate of the states which is used as the initial states, $\hat x_0$, when solving the following OCP for the manipulated inputs.
\begin{subequations}\label{eq:ocp:ms}
	\begin{align}\label{eq:ocp:ms:phi}
		\min_{[x(t)]_{t_0}^{t_f}, \{u_k\}_{k=0}^{N-1}} \quad & \phi = \phi([x(t)]_{t_0}^{t_f}, \{u_k\}_{k=0}^{N-1}),
	\end{align}
	subject to
	\begin{align}
		\label{eq:ocp:ms:x0}
		x(t_0) &= \hat x_0, \\
		\label{eq:ocp:ms:dx}
		\dot x(t) &= f(t, x(t), u(t), d(t), \theta), & t &\in [t_0, t_f], \\
		\label{eq:ocp:ms:u}
		u(t) &= u_k, \quad t \in [t_k, t_{k+1}[, & k &= 0, \ldots, N-1, \\
		\label{eq:ocp:ms:d}
		d(t) &= \hat d_k, \quad t \in [t_k, t_{k+1}[, & k &= 0, \ldots, N-1, \\
		\label{eq:ocp:ms:bounds}
		u_{\min} &\leq u_k \leq u_{\max}, & k &= 0, \ldots, N-1.
	\end{align}
\end{subequations}
The prediction and control horizon, $[t_0, t_f]$, is 6~h, and each of the $N$ control intervals is 5~min. The objective function in~\eqref{eq:ocp:ms:phi} is described in Section~\ref{sec:OptimalControlProblem:Objective:Function}, \eqref{eq:ocp:ms:x0} is the initial condition, \eqref{eq:ocp:ms:dx} is the MVP model where the process noise is disregarded, \eqref{eq:ocp:ms:u}--\eqref{eq:ocp:ms:d} are zero-order-hold parametrizations of the manipulated inputs and the estimated disturbance variables, and \eqref{eq:ocp:ms:bounds} are bounds on the manipulated inputs. Only the first set of manipulated inputs, $u_0$, are administered before a new CGM measurement is received and the horizon $[t_0, t_f]$ is shifted by one control interval.

\subsection{Objective function}\label{sec:OptimalControlProblem:Objective:Function}
The objective function depends on whether insulin or glucagon is administered. In both cases, it is in the form
\begin{equation}\label{eq:ocp:ms:obj}
	\phi = \int_{t_0}^{t_f} \rho_z(z(t)) dt + \sum_{k=0}^{N-1} \rho_{u}(u_k),
\end{equation}
where $\rho_z$ and $\rho_{u}$ are penalty functions and the outputs (the CGM measurements of the blood glucose concentration) are $z(t) = g(t, x(t), \theta)$.

\begin{figure}
	\centering
	\includegraphics[trim=55 2 90 30, clip,width=0.95\linewidth]{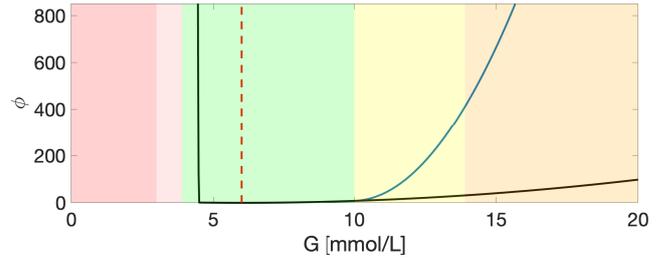}
	\caption{Blood glucose penalty function when administering insulin (blue) and glucagon (black).}
	\label{fig:asymetricQuadraticObjectiveFunctionPlot}
\end{figure}

The penalty function in the first term of~\eqref{eq:ocp:ms:obj} is
\begin{equation}\label{eq:ocp:ms:rhoz}
	\rho_z(z) = \alpha_{\bar{z}} \rho_{\bar{z}}(z) + \alpha_{z_{\min}}\rho_{z_{\min}}(z) + \alpha_{z_{\max}}\rho_{z_{\max}}(z),
\end{equation}
where 1)~the first term penalizes the deviation of the blood glucose from the setpoint $\bar z=6$~mmol/L, 2)~the second term penalizes hypoglycemia ($z < z_{\min}=4.5$~mmol/L), and 3)~the third term penalizes hyperglycemia ($z > z_{\max}=10$~mmol/L):
\begin{subequations}\label{eq:ocp:ms:rhoz:terms}
	\begin{align}
		\label{eq:ocp:ms:rhoz:terms:bar}
		\rho_{\bar{z}} (z) &= \frac{1}{2}(z - \bar{z})^2, \\
		\label{eq:ocp:ms:rhoz:terms:min}
		\rho_{z_{\min}}(z) &= \frac{1}{2}(\min\{0, z - z_{\min}\})^2, \\
		\label{eq:ocp:ms:rhoz:terms:max}
		\rho_{z_{\max}}(z) &= \frac{1}{2}(\max\{0, z - z_{\max}\})^2.
	\end{align}
\end{subequations}
The weights in~\eqref{eq:ocp:ms:rhoz} are $\alpha_{\bar{z}} = 1$, $\alpha_{z_{\min}} = 10^6$, $\alpha_{z_{\max}} = 50$ when computing the insulin flow rates, and $\alpha_{z_{\max}} = 0$ when computing the glucagon flow rate. The penalty function is shown in Fig.~\ref{fig:asymetricQuadraticObjectiveFunctionPlot}. As is evident, preventing hypoglycemia has the highest priority.

When computing the insulin flow rates, the penalty function in the second term of~\eqref{eq:ocp:ms:obj} is
\begin{align}
	\rho_u(u_k) = \rho_{u, ba}(u_{ba, k}) + \rho_{u, bo}(u_{bo, k}),
\end{align}
where
\begin{subequations}\label{eq:ocp:ms:rhou}
	\begin{align}
		\label{eq:ocp:ms:rhou:ba}
		\rho_{u,ba}(u_{ba,k}) &= \lVert u_{ba,k} - \bar{u}_{ba,k} \rVert_2^2, \\
		\label{eq:ocp:ms:rhou:bo}
		\rho_{u,bo}(u_{bo,k}) &= \lVert u_{bo,k} \rVert_1,
	\end{align}
\end{subequations}
in order to penalize excursions from the nominal basal rate $\bar u_{ba, k}$ and to promote the administration of fewer larger insulin boluses.
When computing the glucagon flow rate,
\begin{align}
	\rho_u(u_k) = \lVert u_{G, k} \rVert_2^2,
\end{align}
in order to minimize the administered glucagon.

\subsection{Numerical solution}\label{sec:OptimalControlProblem:Numerical:Solution}
We use a multiple-shooting approach~\citep{Bock:Plitt:1984} to transcribe~\eqref{eq:ocp:ms} by discretizing the dynamic constraint~\eqref{eq:ocp:ms:dx} and the integral in~\eqref{eq:ocp:ms:obj} using an explicit Runge-Kutta method with fixed step size. The result is a nonlinear program which we solve using a sequential quadratic programming~(SQP) method~\citep{Nocedal:Wright:2006}.

\section{Heuristics}
\label{sec:Heuristics}
In the NMPC algorithm, we use a set of heuristics to 1)~switch between administering insulin and glucagon, 2)~compute upper bounds on the manipulated inputs, 3)~modify algorithmic hyperparameters during exercise, and 4) post-process the manipulated inputs computed by solving the OCP~\eqref{eq:ocp:ms}.

\subsection{Switching between insulin and glucagon administration}
Insulin and glucagon should not be administered simultaneously. Therefore, we switch between them based on the blood glucose concentration. If the blood glucose concentration becomes lower than 4.5~mmol/L, the AP switches to glucagon administration. Conversely, if it becomes higher than 5~mmol/L, the AP switches to administering insulin. However, for 1~h after each meal, only insulin can be administered.

\subsection{Bounds and meal-specific heuristics}
The upper bounds on the insulin and glucagon boli are updated at the beginning of every control interval (i.e., whenever a CGM measurement is obtained). Furthermore, the upper bound on the insulin basal rate is twice the target basal rate, i.e., $2 \bar u_{ba, k}$, and the lower bounds on all three manipulated inputs are 0. 

The upper bound on the insulin bolus is
\begin{equation}
	u_{bo, k}^{\max} = \max\{\epsilon, u_{bo, k}^\mathrm{corr} + u_{bo, k}^\mathrm{meal} - u_{bo, k}^\mathrm{hist}\},
\end{equation}
where $\epsilon = 10^{-3}$~[--], $u_{bo, k}^\mathrm{corr}$ is the maximum correction bolus infusion rate, $u_{bo, k}^\mathrm{meal}$ is the maximum meal bolus infusion rate, and $u_{bo, k}^\mathrm{hist}$ is the sum of the insulin bolus infusion rates administered during the previous 11 control intervals. If a meal was announced at time $t_k$, or if no meal was announced in the last hour,
\begin{equation}
	u_{bo, k}^\mathrm{corr}
	= \max\bigg\{0, \frac{1}{T_s}\frac{G - 10~\text{mmol/L}}{ISF}\bigg\},
\end{equation}
where $ISF$~[(mmol/L)/mU] is the insulin sensitivity factor and $T_s = 5$~min is the sampling time. Otherwise, $u_{bo, k}^\mathrm{corr}~=~u_{bo, k-1}^\mathrm{corr}$.
If a meal was consumed within the last hour, the maximum meal bolus is
\begin{align}
	u_{bo, k}^\mathrm{meal}
	&= \max\bigg\{0, \frac{\gamma}{T_s}\frac{\hat{d}}{ICR}\bigg\},
\end{align}
where $\hat d$~[g CHO] is the announced carbohydrate content of the last meal, $\gamma=1.15$~[--] is a bolus allowance factor, and $ICR$~[g/mU] is the insulin-to-carb ratio. Otherwise, $u_{bo, k}^\mathrm{meal} = 0$.
Finally, the insulin bolus history is
\begin{align}
	u_{bo, k}^\mathrm{hist} &= \sum_{j=1}^{11} u_{bo, k-j|k-j},
\end{align}
where $u_{bo, k|k}$ is the insulin bolus infusion rate in the $k$'th control interval. When a meal is announced, we set $u_{bo, k}^\mathrm{hist} = 0$. The history spans 11 control intervals in order to bound the amount of bolus insulin over each 1~h period.

The maximum glucagon bolus is computed by
\begin{equation}
	u_{G, k}^{\max} = \max\{\epsilon, \bar u_G^{\max} - u_{G, k}^\mathrm{hist}\},
\end{equation}
where $\bar u_G^{\max} = 300~\mu$g and the glucagon bolus history is
\begin{align}
	u_{G, k}^\mathrm{hist} &= \sum_{j=1}^{23} u_{G, k-j|k-j}.
\end{align}
Here, $u_{G, k|k}$ is the glucagon infusion rate in the $k$'th control interval. As for the insulin bolus history, the glucagon history spans 23 control intervals in order to bound the amount of glucagon administered over each 2~h period.

Finally, to avoid that the insulin sensitivity is adjusted after a meal, the insulin sensitivity diffusion coefficient is set to zero if a meal was consumed within the last hour. Otherwise, it is set to the value estimated during the parameter estimation. Additionally, in the CD-EKF, we set the state variance of $\log S_I$ as well as the corresponding covariances with the other states to zero when a meal is announced, and we enforce that
\begin{align}
	\log S_I(0) - 1 \leq \log S_I(t) \leq \log S_I(0) + 1
\end{align}
using clipping, where $\log S_I(0)$ is estimated during the parameter estimation.

\subsection{Exercise logic}
During physical activity, the setpoint, $\bar z$, is increased from 6~mmol/L to $7$~mmol/L, the glucagon switching threshold is increased to $7$~mmol/L, and a glucagon dose of $100$~$\mu$g is administered if the blood glucose concentration is below $7$~mmol/L when the physical activity is initiated.

\subsection{Post-processing and open-loop fallback strategy}
Once the solution to the OCP~\eqref{eq:ocp:ms} has been obtained, the resulting manipulated inputs in the first control interval are rounded to the pump resolution which is 0.01~U/h for the insulin basal rate, 0.1~U for the bolus insulin, and 0.01~$\mu$g/h for the glucagon infusion rate.

Furthermore, the AP algorithm is intended to be used by real people, e.g., in clinical trials. Therefore, as a safety measure, we implement the following open-loop strategy if unforeseen circumstances prevent the solution of the OCP~\eqref{eq:ocp:ms}.
%
\begin{subequations}\label{eq:open:loop:fallback}
	\begin{align}
		\label{eq:open:loop:fallback:ba}
		u_{ba, k|k} &=
		\begin{cases}
			0 										& G \leq 8.0~\text{mmol/L}, \\
			\bar{u}_{ba, k} & \text{otherwise},
		\end{cases} \\
		\label{eq:open:loop:fallback:bo}
		u_{bo, k|k} &= 0, \\
		\label{eq:open:loop:fallback:G}
		u_{G, k|k}  &=
		\begin{cases}
			\min\{15~\mu\text{g}, u_{G, k}^{\max}\} 	& G < 4.5~\text{mmol/L}, \\
			0				& \text{otherwise}.
		\end{cases}
	\end{align}
\end{subequations}
%

\section{Results}
\label{sec:Results}
\begin{figure}
	\centering
	\includegraphics[trim=20 25 50 25, clip,width=1.0\linewidth]{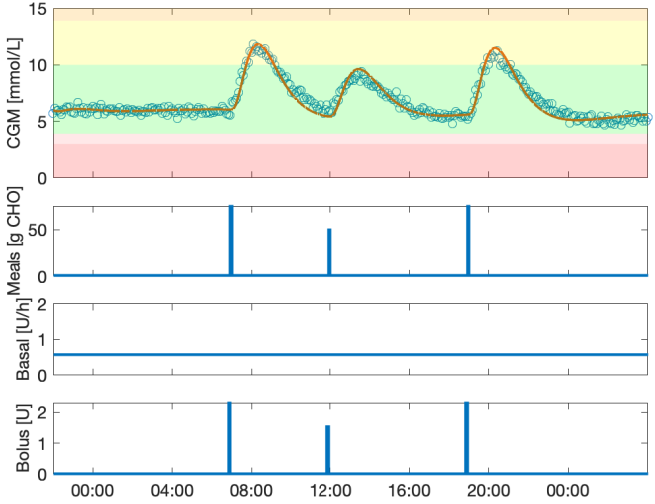}
	\caption{Generated data and simulation with the estimated model for one virtual person. From the top: 1) the CGM data (blue circles) and simulation (red line), 2) the meal carbohydrate content, 3) the insulin basal rate, and 4) the insulin boli.}
	\label{fig:MLEResultsDYCOPS20222ParameterEstimateResults211215f1}
\end{figure}

\begin{figure}
	\centering
	\includegraphics[trim=0 10 0 0, clip,width=1.0\linewidth]{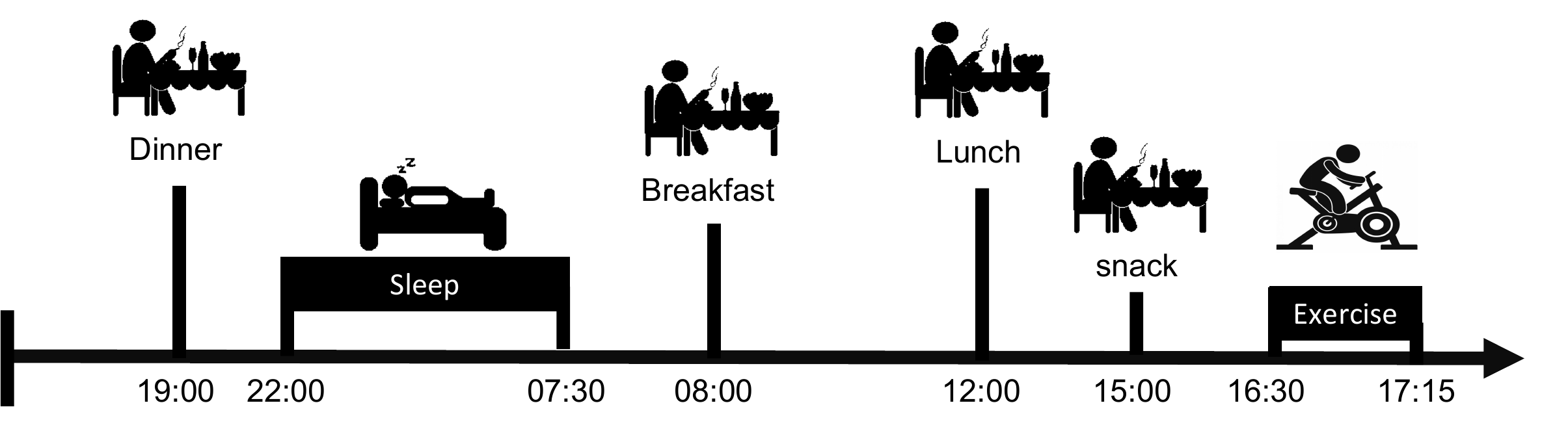}
	\caption{The protocol used in the virtual clinical trial. The protocol consists of a dinner of $75$~g CHO, sleep, a breakfast of $50$~g CHO, a lunch of $75$~g CHO, a snack of $15$~g CHO, and finally, exercise of moderate intensity.}
	\label{fig:AdolescenceStudy}
\end{figure}

\begin{figure}
	\centering
	\includegraphics[trim=38 80 70 50, clip,width=1.0\linewidth]{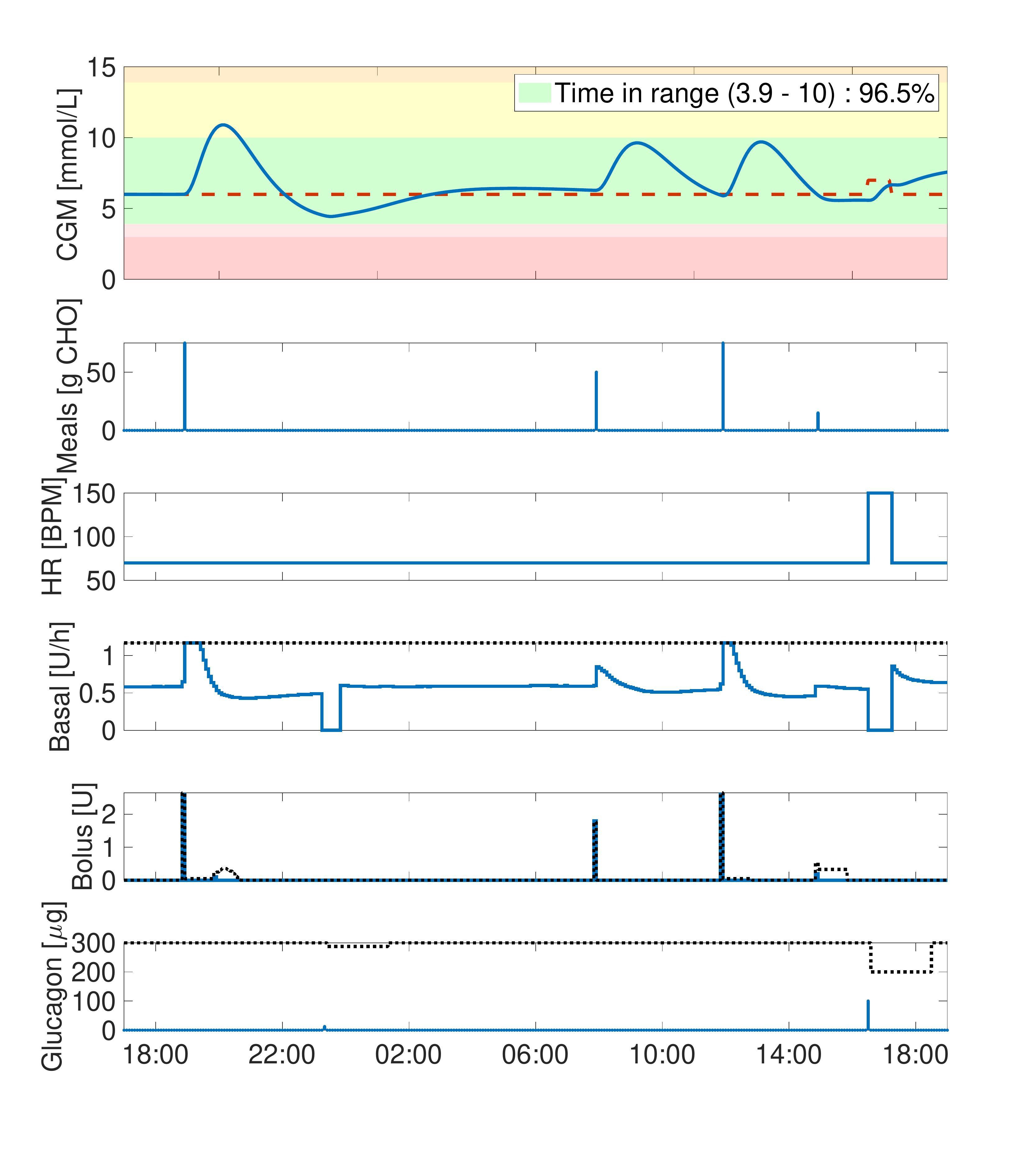}
	\caption{Closed-loop simulation for one virtual person. From the top: 1) The CGM measurement and the setpoint, 2) the meal carbohydrate content, 3) the heart rate, 4) the administered insulin basal rate and the maximum allowed basal rate, 5) the insulin boli and maximum allowed insulin boli, and 6) the glucagon boli and the maximum allowed glucagon boli. The black dotted lines are the upper bounds.}
	\label{fig:AdolescentStudyPatient2DYCOPS2022f1}
\end{figure}

\begin{figure}
	\centering
	\includegraphics[trim=0 30 53 110, clip,width=0.295\linewidth]{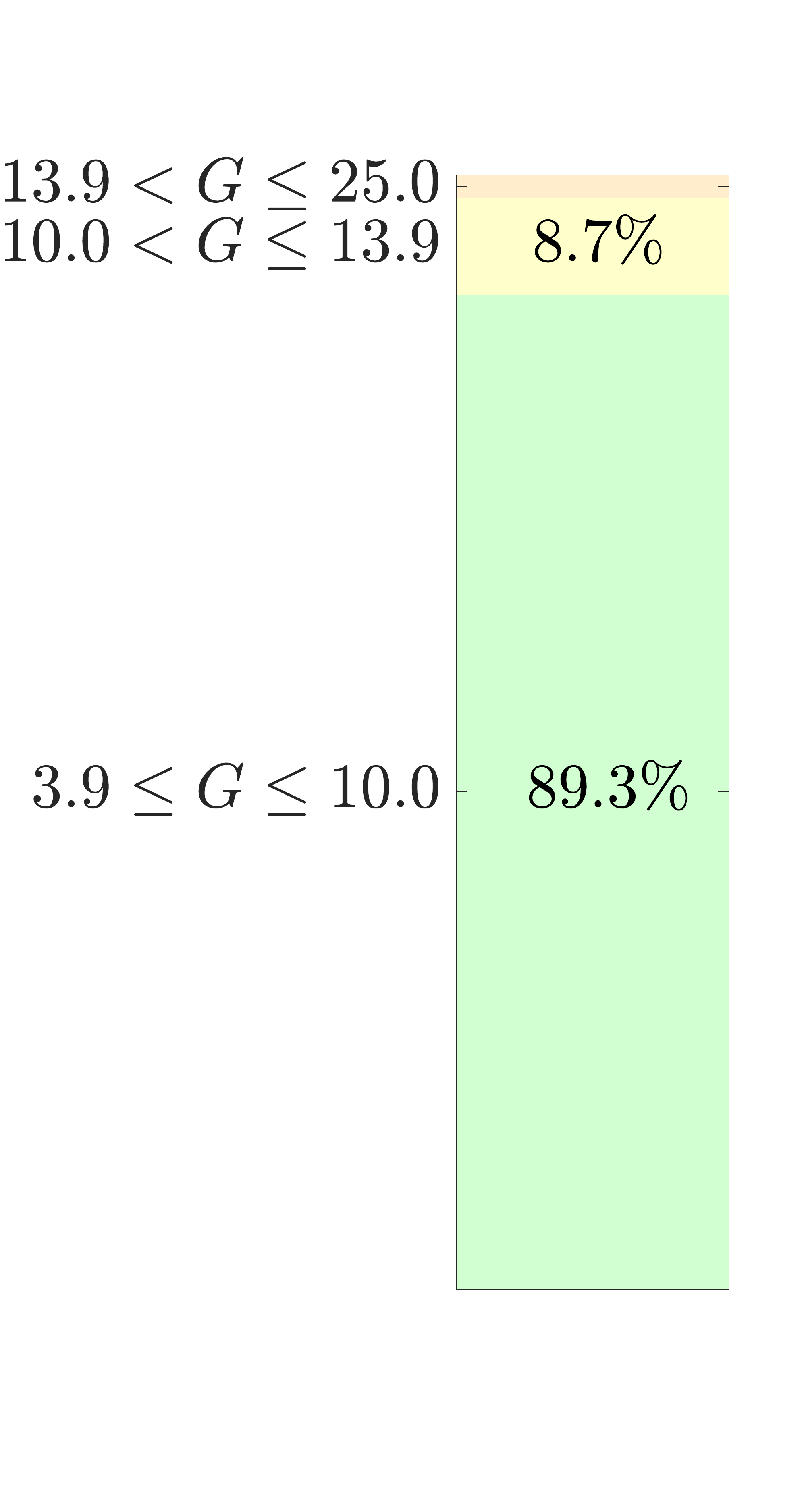}~
	\includegraphics[trim=50 0 60 50, clip,width=0.68\linewidth]{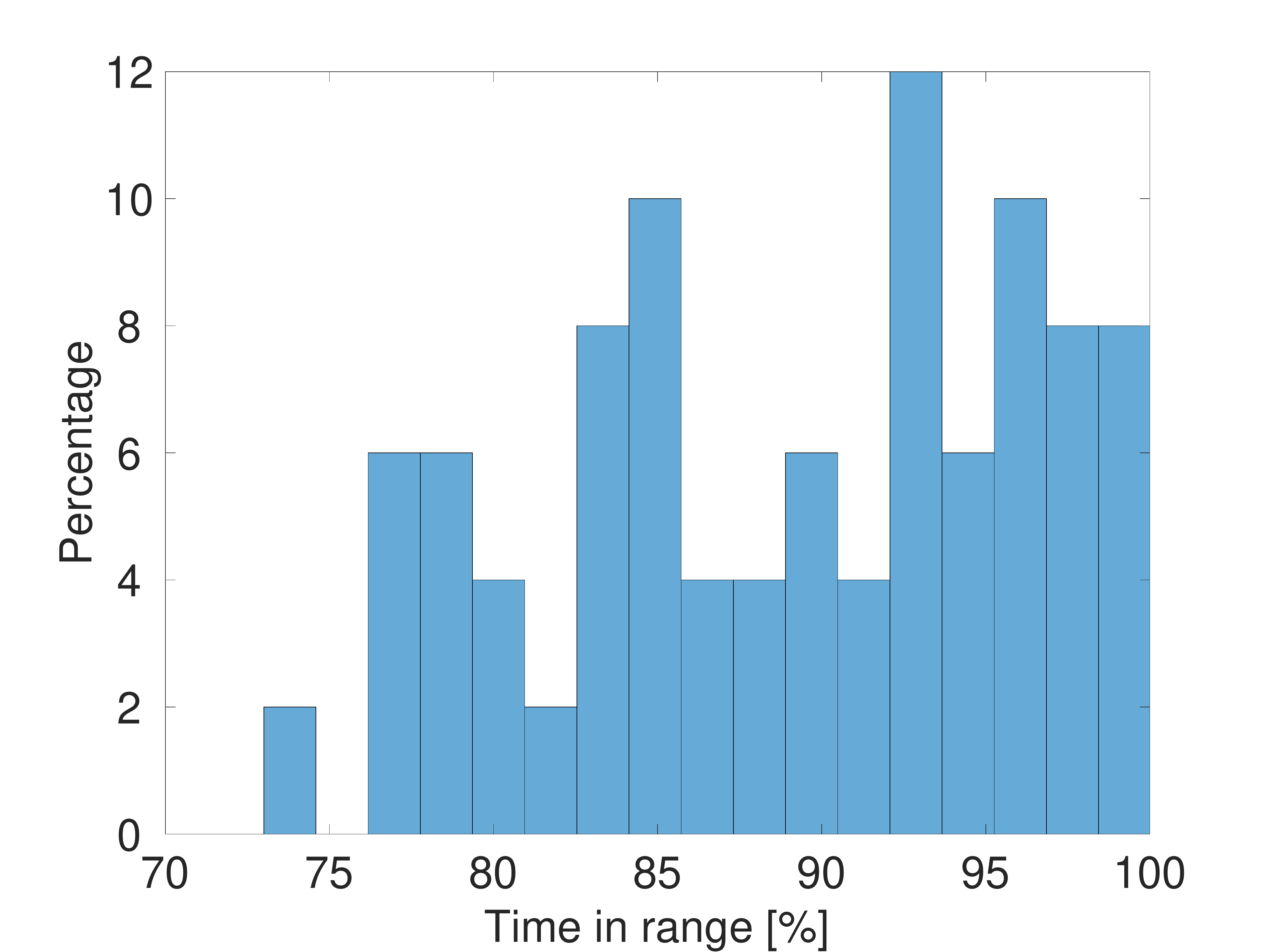}
	\caption{Time in ranges for the 50 virtual people with T1D. Left: mean time in ranges, Right: distribution of time in range.}
	\label{fig:AdolescentStudyAllPatientsDYCOPS2022f3}
\end{figure}

\begin{figure}
	\centering
	\includegraphics[trim=40 52 90 28, clip,width=0.97\linewidth]{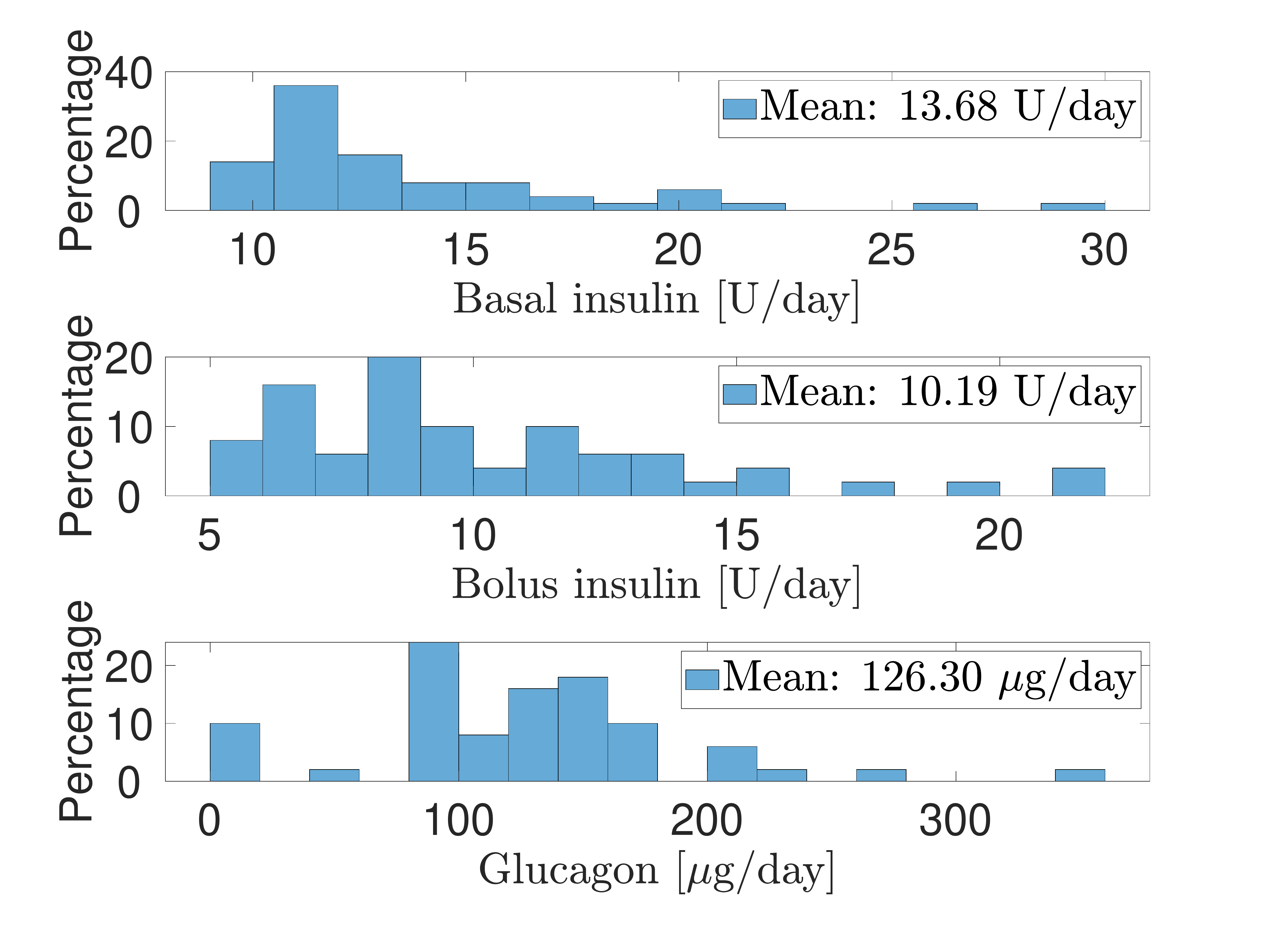}
	\caption{Distributions of the total daily basal and bolus insulin and bolus glucagon for the 50 virtual people. Top: Basal insulin. Middle: Bolus insulin. Bottom: Bolus glucagon.}
	\label{fig:AdolescentStudyAllPatientsDYCOPS2022f4}
\end{figure}

In this section, we present the results from testing both the system identification and the artificial pancreas in a virtual clinical trial with 50 virtual people with T1D. We estimate the parameters in the control model from data generated with the simulation model individually for each person. We estimate the parameters $k_m$, $\tau_D$, $V_G$, $EGP$, $\sigma_G$, and $\sigma_{S_I}$ as well as the initial states in the MVP model. The remaining parameters are fixed. We estimate $ICR$ as described by \cite{Sejersen:etal:2021} and $ISF = 2$~mmol/L/U for all participants. We show the generated data and a simulation (without process noise) with the estimated model in Fig. \ref{fig:MLEResultsDYCOPS20222ParameterEstimateResults211215f1} for one virtual person. In the virtual clinical trial, we use the protocol shown in Fig. \ref{fig:AdolescenceStudy}. In Fig. \ref{fig:AdolescentStudyPatient2DYCOPS2022f1}, we show the closed-loop simulation for the person identified in Fig. \ref{fig:MLEResultsDYCOPS20222ParameterEstimateResults211215f1}. We divide the blood glucose concentration into the following 5 ranges \citep{Holt:etal:2021b} given in mmol/L. Red: severe hypoglycemia (below 3). Light red: hypoglycemia (3--3.9). Green: normoglycemia (3.9--10). Yellow: hyperglycemia (10--13.9). Orange: severe hyperglycemia (above 13.9). For this person, the administered meal bolus is at the limit for most of the meals, but for the snack, we see that only a part of the allowed bolus is administered. The AP is allowed to give the remaining bolus insulin for 1 hour after the snack. The basal rate is increased after the meals where the bolus size is constrained, as the control model predicts the meal to have a larger effect than what the allowed bolus can correct for. Due to the high insulin dose after the dinner, the AP administers a small dose of glucagon. Furthermore, we see that the AP is allowed to give a small correction bolus after the dinner, but once the blood glucose concentration is below $10$~mmol/L, the correction bolus is no longer allowed. It is desired to administer a larger meal bolus and compensate by decreasing the basal rate to reduce the postprandial peak. However, for safety reasons, we have a limit on the maximum bolus even though it can decrease the performance of the AP for some people. After the person begins to exercise, the setpoint is increased and a glucagon bolus of $100$~$\mu$g is administered.\\
Fig. \ref{fig:AdolescentStudyAllPatientsDYCOPS2022f3} and \ref{fig:AdolescentStudyAllPatientsDYCOPS2022f4} show the TIR for all 50 virtual people. The average time in normoglycemia is very high at $89.3$\%, with $8.7$\% in hyperglycemia, $2$\% in severe hyperglycemia and no time in hypoglycemia. The person with the lowest time in normoglycemia still spends more than $70$\% of the time in normoglycemia which is the minimum recommended by \cite{Holt:etal:2021b}. From Fig. \ref{fig:AdolescentStudyAllPatientsDYCOPS2022f4}, we see that the 50 virtual people receive relatively low total daily insulin doses and are sensitive to insulin. The amount of administered glucagon is fairly low and at a reasonable level.

\section{Conclusion}
\label{sec:Conclusion}
In this paper, we present a dual-hormone AP algorithm for controlling the blood glucose concentration in people with T1D. The AP is based on a switching NMPC algorithm, and we use an extension of the MVP model for prediction. Furthermore, we use MLE to estimate the model parameters. The CD-EKF is used in both the NMPC and the MLE algorithm.
We test the algorithm using an extension of the model by~\citet{Hovorka:etal:2002}, and we demonstrate that the average TIR for 50 virtual people with T1D in a virtual clinical trial is 89.3\%. Furthermore, none of the virtual people experience any hypoglycemic events.


\bibliography{bib/ifacconf,ref/References}     

\appendix

\end{document}